\documentclass[journal,twoside,web]{ieeecolor}
\usepackage{generic}

\usepackage{cite}\usepackage{hyperref}
\usepackage{amsmath,amssymb,amsfonts}
\usepackage{algorithmic}
\usepackage{graphicx}
\usepackage{textcomp}
\usepackage{bm,bbm,mathrsfs,amscd}
\usepackage{calc}
\usepackage{color}
\usepackage{dsfont}
\usepackage{graphicx}
\usepackage{epstopdf}
\usepackage{epsfig}
\usepackage{tikz}
\usepackage{amsmath}
\usepackage{amsfonts}
\usepackage{amssymb}
\usepackage{rotating}
\usepackage{mathtools}
\usepackage{color}
\usepackage{subfig}
\usepackage{enumerate,pgfplots}
\newtheorem{theorem}{Theorem}

\newtheorem{proposition}{Proposition}
\newtheorem{lemma}{Lemma}

\newtheorem{example}{Example}
\newtheorem{remark}{Remark}
\newtheorem{assumption}{Assumption}
\newcommand{\ba}{\begin{array}}
\newcommand{\ea}{\end{array}}

\newcommand{\mb}{\boldsymbol}
\newcommand{\be}{\begin{equation}}
\newcommand{\ee}{\end{equation}}

\newcommand{\ds}{\displaystyle}

\newcommand{\eps}{\varepsilon}

\newcommand{\mc}{\mathcal}

\newcommand{\ov}{\overline}

\def\1{\boldsymbol{1}}
\def\0{\boldsymbol{0}}

\newcommand{\R}{\mathbb{R}}

\newcommand{\de}{\mathrm{d}}

\newcommand{\tcb}{\textcolor{black}}

\DeclareMathOperator{\sgn}{sgn}

\def\R{\mathbb{R}}

\tikzstyle{v_c}=[circle, draw,inner sep=2pt, minimum width=12pt, color=blue]
\tikzstyle{v_a}=[circle, draw,inner sep=2pt, minimum width=12pt, color=red]
\tikzstyle{edge} = [draw,thick,-,font=\small ]
\tikzstyle{label} = [draw,fill=black,font=\normalsize]

\def\BibTeX{{\rm B\kern-.05em{\sc i\kern-.025em b}\kern-.08em
	T\kern-.1667em\lower.7ex\hbox{E}\kern-.125emX}}
\begin{document}
\title{On the dynamic behavior of the network SIR epidemic model}
\author{Martina Alutto, Leonardo Cianfanelli, Giacomo Como, \IEEEmembership{Member, IEEE}, Fabio Fagnani
\thanks{The authors are with the  Department of Mathematical Sciences ``G.L.~Lagrange,'' Politecnico di Torino, 10129 Torino, Italy.  E-mail: {\{\!martina.alutto;leonardo.cianfanelli;giacomo.como,fabio.fagnani\!\}@polito.it}. G.~Como is also with the Department of Automatic Control, Lund University, 22100 Lund, Sweden.}
\thanks{This work was partially supported by a MIUR  Research Project PRIN 2017 ``Advanced Network Control of Future Smart Grids'' (http://vectors.dieti.unina.it) and by the Compagnia di San Paolo. }}

\maketitle

\begin{abstract}
We study a susceptible-infected-recovered (SIR) epidemic model on a network of $n$ interacting subpopulations. We analyze the \textcolor{black}{transient and asymptotic} behavior of the \tcb{infection dynamics} in each node of the network. In contrast to the classical scalar epidemic SIR model, where the \tcb{infection curve} is known to \tcb{be unimodal} (either always decreasing over time, or initially increasing until reaching a peak and from then on monotonically decreasing and asymptotically vanishing), we show the possible occurrence of multimodal infection curves in the network SIR epidemic model with $n\ge2$ subpopulations. We then focus on the special case of rank-$1$ interaction matrices, modeling subpopulations of homogeneously mixing individuals with different \tcb{activity rates, susceptibility to the disease, and infectivity levels}. \textcolor{black}{For this special case, we find $n$ invariants of motion and provide an explicit expression for the limit equilibrium point. We also determine necessary and sufficient conditions for stability of the equilibrium points. We then establish }an upper bound on the number of changes of monotonicity of the \tcb{infection curve} at the single node level and provide sufficient conditions \tcb{for its multimodality}. 
Finally, we present some numerical results revealing that, in the case of interaction matrices with rank larger than $1$, the single nodes' infection curves may display multiple peaks.
\end{abstract}

\begin{IEEEkeywords}
Network epidemic models, Susceptible-Infected-Recovered model, infection curves, invariants of motion, limit equilibrium points, stability.
\end{IEEEkeywords}

\section{Introduction}\label{sec:introduction}
The emergence of the COVID-19 pandemic has renewed a huge interest on mathematical models of epidemics. These have proven to be effective tools both for forecasting the spread of infection in the population and for supporting the design of containment rules such as social distantiation and lockdown policies. One of the simplest and most studied of these models is the susceptible-infected-recovered (SIR) epidemic model introduced almost one century ago \cite{Kermack.McKendrick:1927,Hethcote2000TheMO, MKendrickApplicationsOM, Diekmann2000MathematicalEO}. 

In the classical SIR epidemic model, a population is split into three compartments: the \textit{susceptible} individuals, who have not been infected yet and can still catch the disease, the \textit{infected} individuals, who are currently carrying the pathogen and may transmit the disease, and the \textit{recovered} individuals, who have healed from the infection and are forever immune. The rate of new infections is assumed to be proportional to the product between the \tcb{mass} of the susceptible individuals and that of the infected individuals, due to pairwise interactions between them. A crucial index in the analysis of the classical SIR epidemic model is the \emph{reproduction number} $R(t)$, a time-dependent scalar quantity describing the average number of new infections that an infected individual is currently causing. \tcb{Specifically, the fraction of infected individuals at time $t$ is decreasing when $R(t)<1$, while it is increasing for $R(t)>1$.}
 As the reproduction number can be proved to be monotonically decreasing in time and to eventually achieve values below $1$, the infection curve in the classical SIR epidemic model is necessarily \emph{unimodal}. Precisely, if $R(0)\le1$, then the fraction of infected individuals is always monotonically decreasing and vanishes asymptotically. On the other hand, if $R(0)>1$, then the infection curve is initially increasing up to reaching a \emph{peak} value at the time $\hat t$ when $R(\hat t\,)=1$ and is monotonically decreasing from then on, again vanishing as time gets large. This dichotomy has been shown to hold true also for generalizations of the classical SIR epidemic model accounting for more complex interaction mechanisms \cite{CAPASSO197843,Alutto2021OnSE}. It is also at the basis of several control strategies, including some recently proposed in the context of the COVID-19 pandemic, see, e.g., \cite{Cianfanelli.ea:2021} and  \cite{Miclo.ea:2022}. 

The classical SIR epidemic model relies on a number of homogeneity assumptions on the population regarding the  \tcb{individuals' mixing and contact frequency, their aptitude to contract and spread the infection, as well as the time needed to recover. As these assumptions can hardly be met in realistic scenarios, this has motivated the introduction of network versions of the epidemic models \cite{ pastor2015epidemic, pare2020modeling, Zino.Cao:2021}. In this framework, the nodes of the network represent either single individuals or subpopulations of indistinguishable individuals. In the latter case, different subpopulations correspond to, e.g., different geographical areas or age groups \cite{hethcote1978immunization, Fall:2007, Nowzari.ea:2016, Mei.ea:2017, Ogura2016StabilityOS}.}
\tcb{In this paper, we study network SIR epidemic models whereby the nodes of the network are identified with $n$ subpopulations of homogeneous individuals.
{The features of these subpopulations are} encoded into an \emph{interaction matrix} $A$, whose entries $A_{ij}$ represent the rate of new infections in subpopulation $i$ due to the presence of infected individuals in subpopulation $j$ and may incorporate the peculiar susceptibility of individuals in $i$, the infectivity of individuals in $j$, the size of subpopulation $j$ and the rate of interactions among members of the two subpopulations.} 

\tcb{ Most of the literature on network epidemic models focuses on variants of the so called susceptible-infected-susceptible (SIS) epidemic model, differing from the SIR epidemic model in that recovered individuals are again susceptible. In the seminal paper \cite{lajmanovich1976deterministic}, a fundamental analysis of the network SIS epidemic model was carried on, showing the existence of a bifurcation: when the dominant eigenvalue $\lambda_{\max}(A)$ of the interaction matrix is less than or equal to the recovery rate $\gamma$, the disease-free state is a globally asymptotically stable equilibrium point, whereas when $\lambda_{\max}(A)>\gamma$ there exists an endemic equilibrium point that attracts all initial states except for the disease-free state, which remains an unstable equilibrium point. 
It is worth pointing out that the network SIS epidemic model with $n$ subpopulations is an $n$-dimensional monotone system \cite{HirschMonotoneSystems06}, a key property that is at the basis of the results proved in \cite{lajmanovich1976deterministic}.  
More recently, generalizations of the network SIS epidemic model have been considered accounting for, e.g., high-order interactions \cite{cisnerosvelarde2021multigroup}, coupling with opinion dynamics \cite{Xuan2020OnAN,paarporn2017networked}, or competing viruses \cite{anderson2023equilibria, liu2019analysis}.} 

\tcb{In contrast to the aforementioned network SIS epidemic model, the network SIR epidemic model with $n$ subpopulations is a $2n$-dimensional dynamical system and it is not monotone. Such higher dimensionality and lack of monotonicity imply significant analytical challenges. In fact, despite its numerous applications, the network SIR epidemic model is much less studied in the literature. Recent papers such as \cite{acemoglu2021optimal} use calibrated network SIR epidemic models to examine the impact of age-targeted mitigation policies for the COVID-19 pandemic showing how such policies (even with just two age groups) can outperform uniform intervention policies in terms of both mortality rates and economic productivity. 
While most of the studies on the network SIR epidemic model are empirical, there are two notable exceptions. In \cite{Mei.ea:2017}, the authors find a network reproduction number that is a decreasing function of time and plays a role similar to the classical SIR epidemic model: when this reproduction number is less than or equal to $1$, then a certain aggregate infection index (a linear combination of the fraction of infected individuals in the various subpopulations) decreases, whereas when this number is above $1$, this aggregate infection index will first increase and, once the reproduction number becomes smaller than $1$, it will start decreasing to $0$.  However, this aggregate infection index is defined through weights that depend on the initial state and are possibly time varying and this severely limits its applicability. In \cite{ellison2020implications}, a network SIR epidemic model with symmetric rank-$1$ interaction matrices is studied, assuming that individuals are equally susceptible and contagious, {have different activity rates, and there is no homophily in the society}. It is shown that heterogeneity may make the system  reach herd immunity with an aggregate of infected individuals smaller than in the scalar SIR epidemic model. To the best of our knowledge, no results are available on the stability of equilibria, invariants of motion, or the behavior of the infection curve for the single nodes, which is relevant in understanding the effectiveness of targeted interventions.}

This paper provides novel theoretical contributions to the understanding of the dynamic behavior of the network SIR epidemic model  with $n$ nodes. We focus on rank-$1$ interaction matrices, a relevant special case previously studied in \cite{acemoglu2021optimal, ellison2020implications}. \tcb{Our contribution is four-fold}. First, we individuate a novel weighted aggregate infection index always exhibiting a unimodal behavior as function of time. Second, we determine $n$ invariants of motion \textcolor{black}{that allow us to derive the limit equilibrium point as function of the initial state and to analyze the stability of such equilibrium points.} Third, we \tcb{analyze the transient behavior} at the single node level, proving that the infection curve at every node can undertake at most two changes of monotonicity before the reproduction number gets below $1$, and since then it is monotonically decreasing. Fourth, we exhibit a class of network SIR epidemic models with just two nodes where the infection curve at one of the two nodes effectively presents a bimodal behavior with two peaks.
\textcolor{black}{Some of the results appeared in a preliminary form in \cite{9992408}, where the analysis was restricted to two-nodes networks with homogeneous subpopulations and contained no results on invariants of motion, nor any characterization of the limit equilibrium points and their stability.}


The rest of the paper is organized as follows. In Section \ref{sec:2} we introduce the network SIR epidemic model and summarize some known results. \tcb{Our main results 
are presented in Sections \ref{sec:3} (unimodality in the weighted aggregate infection curve, invariants of motion, limit values, and stability) and \ref{sec:4} (dynamic behavior at the single node level). }
 In Section \ref{sec:5}, we illustrate numerical simulations on more general networks. In Section \ref{sec:conclusion}, we discuss future research lines.

\subsection{Notation}
We briefly gather here some notational conventions adopted throughout the paper. We denote by $\R$ and $\R_{+}$ the sets of real and nonnegative real numbers. 
The all-1 vector and the all-0 vector are denoted by $\1$ and $\0$ respectively. 
The transpose of a matrix $A$ is denoted by $A^T$. 
\tcb{For $x$ in $\R^n$, let $||x||_1=\sum_i|x_i|$ and $||x||_{\infty}=\max_i|x_i|$ be its $l_1$- and $l_{\infty}$- norms}, while $[x]$ denotes the diagonal matrix whose diagonal coincides with $x$. For an irreducible matrix $A$ in $\R_+^{n\times n}$, we let $\lambda_{max}(A)$ and $v_{max}(A)$ denote respectively the dominant eigenvalue of $A$ and the corresponding left eigenvector normalized in such a way that $\1'v_{max}(A)=1$, {which has positive entries and is unique due to the Perron-Frobenius theorem.}
 Inequalities between two vectors $x$ and $y$ in $\R^n$ are meant to hold true entry-wise, i.e., $x \le y$ means that $x_i\le y_i$ for every $i$, whereas $x< y$ means that $x_i< y_i$ for every $i$, and $x\lneq y$ means that $x_i\le  y_i$ for every $i$ and $x_j<y_j$ for some $j$. 

\section{Network SIR epidemic model}\label{sec:2}
In this section, we introduce the network SIR epidemic model and gather some known results that will prove useful in the sequel. 
We model networks as finite weighted directed graphs $\mc G=(\mc V,\mc E,A)$, where $\mc V=\{1,2,\ldots,n\}$ is the set of nodes, 
$\mc E\subseteq\mc V\times\mc V$ is the set of directed links, 
and $A$ in $\R_+^{n\times n}$ is a nonnegative weight matrix, to be referred as the interaction matrix, with the property that $A_{ij}>0$ if and only if  there exists a link $(i,j)$ in $\mc E$ directed from node $i$ to node $j$. A network is connected if its interaction matrix $A$ is irreducible. 

In a network SIR epidemic model, $n$ interacting subpopulations $i=1,\ldots,n$ are identified with the nodes of a network $\mc G=(\mc V,\mc E,A)$. For every subpopulation $i$, the time-varying variables $x_i=x_i(t)$, $y_i=y_i(t)$ and $z_i=z_i(t)$ represent the fractions of susceptible, infected, and recovered individuals, respectively, so that sum $x_i+y_i+z_i=1$ remains constant in time. 
\tcb{
The entries $A_{ij}$ of the interaction matrix {account for} the contact frequency between individuals of subpopulation $i$ and individuals of subpopulation $j$, the susceptibility of subpopulation $i$, and the infectivity and size of subpopulation $j$.} Finally, a positive scalar parameter $\gamma$ models the recovery rate, which is assumed to be homogeneous across the network.  

The  network SIR epidemic model with interaction matrix $A$ in $\R_+^{n\times n}$ and recovery rate $\gamma>0$ is then the autonomous system of ordinary differential equations
\begin{equation} 
	\label{eq:network-SIR}
	\begin{cases}
		\dot{x}_i = - x_i\sum_j A_{ij} y_j\\
		\dot{y}_i =   x_i \sum_j A_{ij} y_j - \gamma y_i
	\end{cases}	\qquad i=1,\ldots,n\,.
\end{equation}
\tcb{We omit the equation for the evolution of the fraction of recovered individuals $z_{i} =1-x_{i} - y_{i}$, as it can be directly determined by \eqref{eq:network-SIR}.}
The network SIR epidemic model \eqref{eq:network-SIR} can then be more compactly rewritten in the following vectorial form
\begin{equation}
	\label{eq:network-SIR-compact}
		\dot{x} = -   [x] A y\,,\qquad 
		\dot{y} =   [x] A y - \gamma y\,,
\end{equation}
where $x$ and $y$ in $\mathbb{R}^{n}_{+}$ denote the vectors of the fraction of susceptible and infected individuals, respectively,  in the different subpopulations. 
The following result gathers some basic properties of the network SIR epidemic model. 

\begin{proposition}\label{lemma:basic}
Consider the network SIR epidemic model \eqref{eq:network-SIR-compact} with irreducible interaction matrix $A$ in $\R_+^{n\times n}$ and recovery rate $\gamma>0$. Then,  
\begin{enumerate}
\item[(i)]  the set $$\mc S=\left\{(x,y)\in[0,1]^{2n}:x+y\le\1\right\}$$ is \tcb{positively} invariant;
\item[(ii)] the set of equilibrium points in $\mc S$ is $$\mc S^*=\left\{(x^*,\0):\,x^*\in[0,1]^{n}\right\}\,;$$
\item[(iii)] an equilibrium point $(x^*,\0)$ in $\mc S^*$ is unstable if 
$\lambda_{\max}([x^*] A)>\gamma\,.$
\end{enumerate}
Moreover, for every initial state $(x(0),y(0))$ in $\mc S$: 
\begin{enumerate}
\item[(iv)] for every $i=1,\ldots,n$,  $x_i(t)$ is non-increasing for $t\ge0$, and $x_i(0)>0$ if and only if $x_i(t)>0$ for every $t\ge0$; 
\item[(v)] if $y(0)\gneq\0$, then $y(t)>\0$ for every $t>0$; 
\item[(vi)] there exists $x^*$ in $\R_+^n$ such that $\0\le x^*\le x(0)$ and \be\label{eq:limit-SIR}\lim_{t\to+\infty}x(t)=x^*\,,\qquad \lim_{t\to+\infty}y(t)=\0\,.\ee
\end{enumerate}
\end{proposition}\smallskip
\begin{proof} See \cite{Mei.ea:2017}. 
\end{proof}\medskip

\tcb{In the special case when $n=1$, so that the interaction matrix reduces to a positive scalar value $A=\beta>0$, 
the network SIR epidemic model  \eqref{eq:network-SIR-compact}  reduces to the classical scalar SIR epidemic model} 
\be\label{scalar-SIR} \dot x=-\beta xy\,,\qquad \dot y=(\beta x-\gamma) y\,.\ee 
For the scalar SIR epidemic model \eqref{scalar-SIR}, further results are available, including the following fundamental one. 
\begin{proposition}\label{prop:scalar-SIR}
Consider the scalar SIR epidemic model \eqref{scalar-SIR}, with $\beta>0$,  $\gamma>0$,  and 
  initial state $(x(0),y(0))$ such that $0<x(0)\le1-y(0)\le1$. 
Then, 
\begin{enumerate}
\item[(i)] the quantity $\beta(x+y)-\gamma\log x$ is an invariant of motion; 
\item[(ii)] $x^*=\lim\limits_{t\to+\infty}x(t)$ is the unique solution of the equation 
\be\label{eq:x*scalar}\beta x^*-\gamma\log x^*=\beta (x(0)+y(0))-\gamma\log x(0)\,,\ee 
in the interval $(0,\gamma/\beta]$. 
\end{enumerate}
Moreover, if $y(0)>0$, then: 
\begin{enumerate}
\item[(iii)] if $\beta x(0)\le\gamma$, then $y(t)$ is strictly decreasing for $t\ge0$; 
\item[(iv)] if $\beta x(0)>\gamma$, then there exists a peak time $\hat t>0$ such that $y(t)$ is strictly increasing for $t$ in $[0,\hat t]$ and strictly decreasing for $t$ in $[\hat t,+\infty)$.
\end{enumerate}
\end{proposition}\smallskip
\begin{proof}See \cite[Chapter~2.4]{Brauer2019MathematicalMI}. \end{proof}\medskip

We now provide a simple example of a network SIR epidemic model with just $n=2$ nodes where the infection curve at the single node level achieves multiple peaks. \tcb{This contrasts the unimodality  of the infection curve in the scalar SIR epidemic model determined by the dichotomy in Proposition \ref{prop:scalar-SIR}(iii)--(iv).}


\begin{example}\label{example:mutimodal}
\begin{figure}
	\centering
	\includegraphics[scale=0.65]{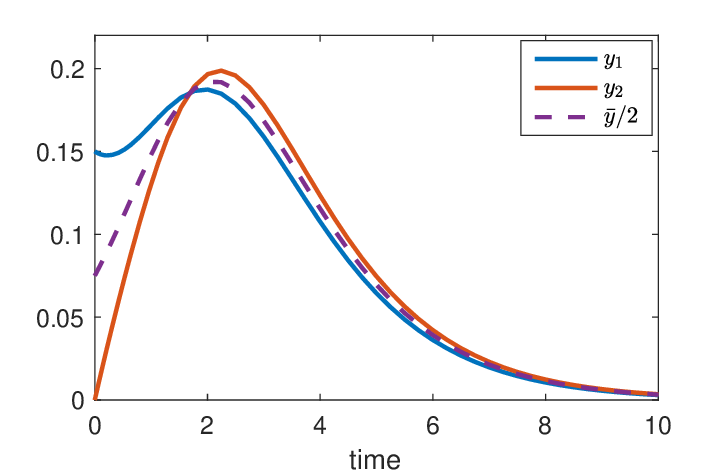}
	\caption[]{Numerical simulation of the network SIR epidemic model with $n=2$ nodes with interaction matrix $A = \1\1^T$, recovery rate $\gamma =1$, and initial state  $y_{1}(0)=1-x_1(0)=\tcb{0.15}$ and $y_{2}(0)=1-x_{2}(0)=0$ satisfying \eqref{eq:conditions}-\eqref{oveps}.}	
	\label{fig:fig1}
\end{figure}
Consider the network SIR epidemic model \eqref{eq:network-SIR} with $n=2$ subpopulations, interaction matrix $A=\1\1^T$, and unitary recovery rate $\gamma=1$.  Let the initial state be 
	\be \label{eq:conditions}y_1(0)=1-x_1(0)=\eps \,,\quad y_2(0)=1-x_2(0)=0\,,\ee
	for some $\eps>0$ such that 
	\be \label{oveps}\frac{1-\eps}{2-\eps}(1-\log(2-\eps))>\eps\,.\ee
The range of such values of $\eps$ is nonempty since the function $$g(\eps)=\frac{1-\eps}{2-\eps}(1-\log(2-\eps))-\eps$$  is continuous in the interval $[0,1]$ and $g(0) = \frac12(1-\log2)>0$.  

Observe that, with these initial states, we have 
\be\label{doty1<0}\dot{y}_{1}(0) = x_{1}(0)(y_1(0)+y_2(0)) - y_{1}(0)=-\eps^2<0\,,\ee
which implies that $y_1(t)$ is strictly decreasing for sufficiently small $t>0$. We now show that $y_1(t)$ cannot remain decreasing for all values of $t\ge0$, but will necessarily become strictly increasing in a certain time range, before eventually starting to decrease again and vanish as $t$ grows large. 

Towards this goal, first observe that the aggregate variables $\ov x=x_1+x_2$ and $\ov y=y_1+y_2$ satisfy an autonomous scalar SIR epidemic model
\be\label{SIR-aggregate}\dot{\ov x}=-\ov x\, \ov y\,,\qquad \dot{\ov y}=(\ov x-1)\ov y\,.\ee
Then, since 
$\dot{\ov y}(0) = (\ov x(0) -1)\ov y(0) >0\,,$
Proposition \ref{prop:scalar-SIR}(iv) implies that there exists a peak time $\hat t>0$ at which $\dot{\Bar{y}}(\hat t\,)=0$, i.e., $\ov{x}(\hat t\,)=1$. 
This fact, Proposition \ref{prop:scalar-SIR}(i), and \eqref{eq:conditions} imply that 
\be\label{eq:bar_y_t*}
\ov y(\hat t\,)
=\ov{x}(0)+\ov{y}(0)-\ov x(\hat t\,)+\log\frac{\ov{x}(\hat t\,)}{\ov{x}(0)}\\
=1-\log(2-\eps)\,.\ee
Since $\dot x_2=-x_2\ov y$ and $x_2(0)=1$, we have that  
\be\label{x2*}x_2(\hat t\,)=\exp\Bigg(-\int_0^{\hat t}\ov y(t)\de t\Bigg)=\frac{\ov x(\hat t\,)}{\ov x(0)}=\frac{1}{2-\eps}\,,\ee
where the second equality follows from integrating the first equation in \eqref{SIR-aggregate} and the last one follows from \eqref{eq:conditions}. 
%
It then follows from  \eqref{eq:bar_y_t*} and \eqref{x2*} that 
\be\label{doty2hat}
 \dot y_2(\hat t\,)  
 =  x_2(\hat t\,) \ov y(\hat t\,) -y_2(\hat t\,) 
= \frac{1-\log(2-\eps)}{2-\eps} - y_2(\hat t\,)\,.
\ee

Now, assume by contradiction that 
$\dot y_1(t)\le 0$ for all $t\ge0$. 
In particular, this would imply that 
$y_1(\hat t\,)\le y_1(0)=\eps$, so that
$$y_2(\hat t\,)=\ov y(\hat t\,)-y_1(\hat t\,)\ge1-\log(2-\eps)- \eps\,.$$ by \eqref{eq:bar_y_t*}. 
Recalling that $\dot{\ov{y}}(\hat t\,)=0$, substituting  the above in the righthand side of \eqref{doty2hat}, and using \eqref{oveps}, we would then get 
$$\dot y_1(\hat t\,)=\dot{\ov{y}}(\hat t\,)-\dot y_2(\hat t\,)\ge\frac{1-\eps}{2-\eps} \left(1-\log(2-\eps)\right) - \eps>0\,,$$
thus contradicting the assumption that $\dot y_1(t)\le 0$ for $t\ge0$. 
It then follows that there must exist some values of time $t\ge0$ such that $\dot y_1(t)>0$. 
Together with \eqref{doty1<0} and the fact $\lim_{t\to+\infty}y_1(t)=0$ by Proposition \ref{lemma:basic}(vi), this implies that  the infection curve $t\mapsto y_1(t)$ is multimodal. 
In fact, the results to be presented  in Section \ref{sec:4} imply that such behavior is necessarily as illustrated in Figure \ref{fig:fig1}, i.e.,  $y_1(t)$  is strictly decreasing in an interval $[0,\check t_1]$, until reaching a positive local minimum point $\check t_1>0$,  it is then strictly increasing in  an interval $[\check t_1,\hat t_1]$ until reaching a second peak at some time $\hat t_1>\check t_1$, and is eventually strictly decreasing for $t\ge\hat t_1$. 

Notice that the network SIR epidemic model considered in this example can  be interpreted as a scalar SIR epidemic model where a single population of individuals has been split into two equally sized subpopulations having distinct initial states. Specifically, all the initially infected individuals belong to the first subpopulation, while the second one initially contains only susceptible individuals. The parameters are chosen such that if the two subpopulations were isolated, the first one would undertake an exponential decrease to a disease-free state. However, because of the presence of the second subpopulation, the infection can further spread and eventually hit back the first subpopulation, making it undergo a second wave of infection with a second peak (the first one being at time $0$).
\end{example}\medskip

\section{The network SIR epidemic model with rank-$1$ interaction  matrices}\label{sec:3}

 In this section, we study the network SIR epidemic model in the special case when the interaction matrix $A$ is irreducible and has rank-$1$, as per the following equivalent assumption. 
 \begin{assumption}\label{ass:rank1}
 The interaction matrix $A$ satisfies 
\be\label{eq:rank1}A=ab^T\,,\ee
for two vectors $a>\0$ and $b>\0$ in $\R^n$. 
 \end{assumption}

For a rank-$1$ interaction matrix $A=ab^T$, 
the network SIR epidemic model's equations \eqref{eq:network-SIR} can be rewritten as
\begin{equation} 
	\label{eq:dynamics}
	\begin{cases}
		\dot{x}_{i} = -a_{i} x_{i} \bar{y}\\
		\dot{y}_{i} = a_{i} x_{i} \bar{y} - \gamma y_{i}
	\end{cases}	\qquad i=1,\ldots,n\,,
\end{equation}
where 
\be\label{barxy}
\bar{y} = \sum_{j=1}^{n} b_{j} y_{j}\,.
\ee
\tcb{is the weighted aggregate of infected individuals that drives the rate of new infections in all subpopulations.}
\begin{remark}\tcb{
The vectors $a$ and $b$ in Assumption \ref{ass:rank1} can be given the following epidemiological interpretation.
On the one hand, the $i$-th entry $a_i$ of vector $a$ accounts for the susceptibility level and the activity rate of subpopulation $i$. 
On the other hand, the $i$-th entry $b_i$ of vector $b$ accounts for the infectivity level, the activity rate, and the size of subpopulation $i$. Assumption \ref{ass:rank1} then boils down to that the contact frequency between individuals of subpopulation $i$ and those of subpopulation $j$ is proportional to the activity rates of the two subpopulations. Therefore, $\bar y$ represents the aggregate of infected individuals weighted by their infectivity and activity rates. This model can capture the effect of targeted containment policies as proposed in \cite{acemoglu2021optimal}. However, it cannot take into account the presence of homophily in the interactions.
} 
Notice that the class of network SIR epidemic models satisfying Assumption \ref{ass:rank1} encompass the ones studied in \cite{ellison2020implications} where authors require that $A$ is both rank-$1$ and symmetric. \end{remark}

\subsection{Unimodality of the weighted aggregate infection curve}



\tcb{In this subsection, we study the dynamics of the weighted aggregate of infected individuals $\bar y$, showing that under Assumption \ref{ass:rank1} this quantity always exhibits a unimodal behavior}. 
For our analysis, it is convenient to introduce two different weighted aggregates of susceptible individuals:
%
\be\label{eq:tildex}
\bar{x} = \sum_{j=1}^{n} b_{j} x_{j}\,, \quad \tilde{x} = \sum_{j=1}^{n}  a_{j} b_{j} x_{j}\,.\ee 
\tcb{The following result describes the dynamics of the weighted aggregates $\bar x$ and $\bar y$ of susceptible and, respectively, infected individuals. In particular, it shows that the logarithmic time derivative of $\bar y$ equals the difference between the weighted aggregate of susceptible individuals $\tilde x$ and the recovery rate $\gamma$}.

\begin{lemma}\label{lemma1}
Consider the rank-$1$ connected network SIR epidemic model \eqref{eq:dynamics}. Then,  
\be\label{dotovx+dotovy}\dot{\ov x}=-\ov y\tilde x\,,\qquad\dot{\ov y}=\ov y\left(\tilde x-\gamma\right)\,.\ee
\end{lemma}\smallskip
\begin{proof} See Appendix \ref{sec:proof-lemma1}. 
\end{proof}\medskip
\tcb{The following is the main result of this subsection. It shows that the weighted aggregate infection curve $t\mapsto\bar y(t)$ is always unimodal and that
there is a dichotomy determining the shape of this curve as a function of the initial value of the weighted aggregate of infected individuals $\tilde x(0)$.}

\begin{theorem}\label{theo:bary}
Consider the rank-$1$ connected network SIR epidemic model \eqref{eq:dynamics}. Let the initial state $(x(0),y(0))$ in $\mc S$ be such that \be\label{initial-cond}\0\lneq x(0)\le\1-y(0)\lneq\1\,.\ee
Then, $\tilde x(t)$ is strictly decreasing for $t\ge0$ and 
\be\label{limtildex}\lim_{t\to+\infty}\tilde x(t)<\gamma\,.\ee  
Moreover: 
\begin{enumerate} 
\item[(i)] if \be\label{subcritical}\tilde x(0)\le\gamma\,,\ee then $\ov y(t)$ is strictly decreasing for $t\ge0$; 
\item[(ii)] if \be\label{supercritical}\tilde x(0)>\gamma\,,\ee then there exists $\hat t>0$ such that $\ov y(t)$ is strictly increasing on $[0,\hat t]$ and strictly decreasing on $[\hat t,+\infty)$.
\end{enumerate}
\end{theorem}\smallskip
\begin{proof}
\tcb{By the way $\tilde x$ is defined in  \eqref{eq:tildex} and using the equations in \eqref{eq:dynamics}, we can express its time derivative as}
%
\be\label{dottilldex}\dot{\tilde x}= \sum_{j=1}^{n}  a_{j} b_{j}\dot x_{j}=-\ov y\sum_{j=1}^{n}  a_{j}^2 b_{j}x_{j}\,.\ee
Now, the rightmost inequality in \eqref{initial-cond} and Proposition \ref{lemma:basic}(v) imply that $y(t)>\0$ for every $t>0$, whereas the leftmost inequality in \eqref{initial-cond} and Proposition \ref{lemma:basic}(iv) imply that there exists some $i$ in $\{1,\ldots,n\}$ such that $x_i(t)>0$ for every $t\ge0$. From \eqref{dottilldex}, and since $a>\0$ and $b>\0$, we get that 
\be\label{dottildex<0}\dot{\tilde x}(t)=-\ov y(t) \sum_{j=1}^{n}  a_{j}^2 b_{j}x_{j}(t)\le-a_i^2b_ix_i(t)\ov y(t)<0\,,\ee
for every $t>0$. Inequality \eqref{dottildex<0} implies that $t\mapsto\tilde x(t)$ is strictly decreasing for $t\ge0$. 
Now, let $$\tilde x(\infty)=\lim_{t\to+\infty}\tilde x(t)\,.$$ If $\tilde x(0)\le\gamma$, then $\tilde x(\infty)<\tilde x(0)\le\gamma$, so that \eqref{limtildex} is satisfied. 
On the other hand, if $\tilde x(0)>\gamma$, assume by contradiction that $\tilde x(\infty)\ge\gamma$. Then $\tilde x(t)>\gamma$ for every $t\ge0$, so that, by \eqref{dotovx+dotovy},
$$\dot{\ov y}(t)=\ov y(t)\left(\tilde x(t)-\gamma\right)\ge0\,,\qquad \forall t\ge0\,.$$
The above would imply that $\ov y(t)$ is nondecreasing, so that 
$$\lim_{t\to+\infty}\ov y(t)\ge\ov y(0)>0\,,$$ 
thus contradicting Proposition \ref{lemma:basic}(vi). Hence, we have $\tilde x(\infty)<\gamma$ also when $\tilde x(0)>\gamma$, thus completing the proof of the first part of the statement. 

(i) If $\tilde x(0)\le\gamma$, by point (i) we have that $\tilde x(t)<\gamma$ for $t>0$. Hence, \eqref{dotovx+dotovy} implies that
$$\dot{\ov y}(t)=\ov y(t)\left(\tilde x(t)-\gamma\right)<0\,\qquad \forall t>0\,,$$
thus showing that $\ov y(t)$ is strictly decreasing for $t\ge0$. 

(ii) If $\tilde x(0)>\gamma$, by point (i), $\tilde x(t)$ is strictly decreasing and 
$$\lim_{t\to+\infty}\tilde x(t)<\gamma\,.$$ Then, there necessarily exits a time $\hat t>0$ such that $\tilde x(t)>\gamma$ for $0\le t<\hat t$, $\tilde x(\hat t)=\gamma$, and $\tilde x(t)<\gamma$ for $t>\hat t$. 
It then follows from \eqref{dotovx+dotovy} that $\ov y(t)$ is strictly increasing for $t$ in $[0,\hat t]$ and strictly decreasing for $t$ in $[\hat t,+\infty)$, thus proving (ii).
\end{proof}\medskip

\begin{remark}
It is proven in \cite{Mei.ea:2017}  that, for general irreducible interaction matrices, if $\lambda_{\max}([x(\tau)] A)<\gamma$ for some $\tau \geq 0$, then the weighted aggregate  $v_{max}(\tau)^{T}y(t)$ is monotonically decreasing to $0$, whereas, if $\lambda_{\max}([x(0)] A)>\gamma$, then for small times $v_{max}(0)^{T}y(t)$ grows exponentially fast. 
Notice that $v_{max}(\tau)^{T}y(t)$ explicitly depends on $x(\tau)$: in particular, $v_{max}(0)^{T}y(t)$ may not be unimodal in general. For rank-$1$ interaction matrices $A=ab^T$, we have $v_{max}(t)=b$ for $t\ge0$, so that $v_{max}(\tau)^{T}y(t)=\ov y(t)$ is unimodal by Theorem \ref{theo:bary}. 
	
\end{remark}

\subsection{Invariants of motion, limit equilibria, and stability}
\tcb{
In this subsection, we study invariants of motion, limit equilibria, and stability for the rank-$1$ network SIR epidemic model \eqref{eq:dynamics}. We start with the following technical result. 
\begin{lemma}\label{lemma:varphi}
Fix two vectors $a>\0$ and $b>\0$ in $\R^n$. 
Then, for every $(x,y)$ in $\mc S\setminus\{(x,\0):\tilde x\ge\gamma\}$, 
\begin{enumerate}
\item[(i)] the equation 
\be\label{eq:barx*} \xi=  \sum_{i=1}^n b_ix_i \exp (a_i(x_i-\ov x-\ov y) / \gamma)\,,\ee
{admits exactly one solution in the interval $[0,\ov x]$ that is denoted by $\xi=\varphi(x,y)$;}
%
%
\item[(ii)] $\varphi(x,{y})=\ov x$ if and only if $y=\0$ or $x=\0$; 
\item[(iii)] $\varphi(x,y)=0$ if and only if $x=\0$;  
\item[(iv)]  $\varphi(x,y)$ is continuous on $\mc S\setminus\{(x,\0):\,\tilde x\ge\gamma\}$.
\end{enumerate} 
\end{lemma}\smallskip
\begin{proof} See Appendix \ref{sec:proof-lemmavarphi}.\end{proof}}\medskip

%
%
%
%
%
%

Our next result  generalizes Proposition \ref{prop:scalar-SIR}(i)--(ii), as it characterizes $n$ invariants of motion \textcolor{black}{as well as the dependence of the limit equilibrium point on the initial state} for the network SIR epidemic model with irreducible rank-$1$ interaction matrix. 
\tcb{\begin{theorem}\label{theo:invariant}
	Consider the rank-$1$ connected network SIR epidemic model \eqref{eq:dynamics}. Then: 
		\begin{enumerate}
		\item[(i)] for every $i=1,\ldots,n$, the quantity 
		\be\label{eq:hi-def}{h_i}=x_i\exp(-a_i(\ov x+\ov y)/\gamma)\ee
		is an invariant of motion;
\item[(ii)] for every initial state $(x(0),y(0))$ in $\mc S$, 
\be\label{eq:limit}\!\!\!\!\!\!\!\!\!\!\!\lim_{t\to+\infty}x(t) = \Phi(x(0),y(0))\,,\ee
where, for every $(x,y)$ in $\mc S$ and $1\le i\le n$, 
$$\Phi_i(x,y)=x_i \exp \left({a_i}\big(\varphi(x,y)-\ov x-\ov y\big)/\gamma\right)\,,$$
{where $\varphi(x,y)$ is defined as in Lemma \ref{lemma:varphi}(i);} 
\item[(iii)] an equilibrium point $(x^*,\0)$ in $\mc S^*$  is stable if  \be\label{stability-cond}\tilde x^* = \sum_{j=1}^{n} a_j b_j x_j^*<\gamma\,.\ee
	\end{enumerate}
\end{theorem}\smallskip}
\begin{proof}\tcb{
	(i) 
	Equations \eqref{eq:dynamics} and \eqref{dotovx+dotovy} imply that  
       $$\ba{rcl}\ds 
	\dot{h_i}
	&=&\ds 
	\left(\dot x_i-a_ix_i\left(\dot{\ov x}+\dot{\ov y}\right)/\gamma\right)\exp(-a_i(\ov x+\ov y)/\gamma)
	\\[2pt]
	&=&\ds \left(-a_ix_i\ov y+ a_i x_i \ov y\right)\exp(-a_i(\ov x+\ov y)/\gamma)	\\[2pt]
	&=&0\,,\ea$$
	thus proving (i).
	}

	\tcb{(ii)
	It follows from Proposition \ref{lemma:basic}(vi) that \eqref{eq:limit-SIR} holds true for some $\0\le x^*\le x(0)$, so that point (i) implies that 
	\be\label{eq:limit-hi}h_i(x(0),y(0))=h_i(x(t),y(t))\stackrel{t\to+\infty}{\longrightarrow}h_i(x^*,0)\,.\ee
	Let 
	$\ov x^*=\sum_ib_ix_i^*$. 
	It then follows from \eqref{eq:hi-def} and \eqref{eq:limit-hi} that 
	$$\ba{rcl}x_i^*\exp(-a_i\ov x^*/\gamma)
	&=&h_i(x^*,0)\\[2pt]
	&=&h_i(x(0),y(0))\\[2pt] 
	&=&x_i(0)\exp(-a_i(\ov x(0)+\ov y(0))/\gamma)\,,\ea$$
        so that 
        $$\lim_{t\to+\infty}x_i(t) =x_i^*= x_i(0)\exp (a_i (\bar x^*-\ov x(0)-\ov y(0))/\gamma)\,.$$
	Multiplying both sides of the above by $b_i$ and summing up over $i=1,\ldots,n$ 
	 shows that $\ov x^*$ is indeed a solution of \eqref{eq:barx*}. 
	 Now, observe that Proposition \ref{lemma:basic}(iv) ensures that $\bar x(t)$ is non-increasing in $t$, so that necessarily $\ov x^*$ belongs to the interval $[0,\ov x(0)]$. By Lemma \ref{lemma:varphi}(i), $\varphi(x(0),y(0))$ is the unique solution of  \eqref{eq:barx*} in $[0,\ov x(0)]$, so that necessarily $\ov x^*=\varphi(x(0),y(0))$, thus completing the proof of point (ii). }
	 
\tcb{(iii)	Let $(x^*,\0)$ be an equilibrium point satisfying \eqref{stability-cond}.  
By  Proposition \ref{lemma:basic}(ii), we have $\Phi(x^*,\0)=x^*$, while Lemma \ref{lemma:varphi}(iv) implies that $\Phi(x,y)-x$ is continuous in the point $(x^*,\0)$. Hence, for every $\eps>0$ there exists $\delta_\eps>0$ such that 
\be\label{continuity}||(x,y)-(x^*,\0)||_{\infty}<\delta_{\eps}\quad\Longrightarrow\quad||\Phi(x, y)-x||_{\infty}<\eps/2\,.\ee
Observe that it is not restrictive to assume that 
\be\label{delta<eps}\delta_{\eps}\le\eps\min\left\{1/2,{b_*}/{||b||_1}\right\}\,,\qquad b_*=\min\nolimits_i b_i\,.\ee 
For $\eps>0$, consider an initial state $(x(0),y(0))$ in $\mc S$ such that 
\be\label{deltaeps}||x(0)-x^*||_{\infty}<\delta_{\eps}\,,\qquad ||y(0)||_{\infty}<\delta_{\eps}\,.\ee
Observe that Proposition \ref{lemma:basic}(iv) and point (ii) imply that 
\be\label{sandwitch}\Phi(x(0),y(0)) \le x(t)\le x(0)\,,\qquad\forall t\ge0\,.\ee
Therefore, 
\be\label{xinfty}\ba{rcl}
\!\!\!||x(t)-x^*||_{\infty}\!\!\!
&\le&\!\!\!||x(t)-x(0)||_{\infty}+||x(0)-x^*||_{\infty} \\[2pt]
&\le&\!\!\!||\Phi(x(0),y(0))- x(0)||_{\infty}+\delta_{\eps}\\[2pt]
&<&\!\!\!\eps/2+\eps/2\\
&=&\!\!\!\eps\,,\ea\ee
where the first inequality follows from from the triangle inequality, the second one from \eqref{sandwitch} and \eqref{deltaeps}, and the third one from  \eqref{continuity} and \eqref{delta<eps}. 
On the other hand, Theorem \ref{theo:bary}(i) implies that $0\le\ov y(t)\le\ov y(0)$ for all $t\ge0$, so that 
\be\label{yinfty}||y(t)||_{\infty}\le\frac{\ov y(t)}{b_*}\le\frac{\ov y(0)}{b_*}\le\frac{||b||_1}{b_*}||y(0)||_{\infty}<\frac{||b||_1}{b_*}\delta_{\eps}\le\eps\,,\ee
the last two inequalities implied by \eqref{deltaeps} and  \eqref{delta<eps}, respectively. Finally, \eqref{xinfty}--\eqref{yinfty} imply that 
$||(x(t),y(t))-(x^*,\0)||_{\infty}<\eps$, thus proving stability of the equilibrium point $(x^*,\0)$. 
}
\end{proof}\medskip
\begin{remark} \tcb{
Invariants of motion have already been used in the literature on scalar epidemic models in order to determine the limit value of the fraction of susceptible individuals as well as the peak value of the fraction of infected individuals: see, e.g., \cite{feng2007final, sontag2023explicit}. For the classical scalar SIR epidemic model and its variant with equal death and birth rates, this has even lead to exact analytical solutions \eqref{scalar-SIR} \cite{Harko2014ExactAS}. 
The novelty of Theorem \ref{theo:invariant}(i) is that it addresses a network SIR epidemic model with $n$ subpopulations, that is a $2n$-dimensional system and determines $n$ invariants of motion for it. Theorem \ref{theo:invariant}(ii) then uses the form of such invariants of motion to characterize the functional dependance of the limit equilibrium point on the initial state, while Theorem \ref{theo:invariant}(iii) builds on this to provide sufficient conditions for stability. }
\end{remark}\medskip
\begin{remark} \tcb{
Theorem \ref{theo:invariant}(iii) complements Proposition \ref{lemma:basic}(iii). Indeed, if the interaction matrix $A$ satisfies Assumption \ref{ass:rank1}, then  $$\lambda_{\max}([x^*]A)=\lambda_{\max}([x^*]ab')=b'x^*=\tilde x^*\,.$$  
Hence, for connected rank-$1$ network SIR epidemic models, an equilibrium point $(x^*,\0)$ is stable if $\tilde x^*<\gamma$ and unstable if $\tilde x^*>\gamma$ (no equilibrium point can be asymptotically stable as they form  a continuum). 
Notice that, while Proposition \ref{lemma:basic}(iii) is an immediate consequence of the linearization in \cite[Theorem 5.2(ii)]{{Mei.ea:2017}}, Theorem \ref{theo:invariant}(iii) is not since such linearization always has $0$ as an eigenvalue with multiplicity $n$.  }
\end{remark}

\section{Dynamic behavior of the infection in the single subpopulations}\label{sec:4}
In this section, we analyze the infection curves of the sigle subpopulations for the connected rank-$1$ network SIR epidemic model.
\subsection{A classification of possible dynamic behaviors}
\tcb{We first define the quantities
\be\label{widef} w_i=\tilde x-\gamma-a_i\ov y\,,\qquad i=1,\ldots,n\,,\ee
measuring how the rate of new infections in each subpopulation is changing in time.}
\tcb{
Indeed, from the second equation in \eqref{eq:dynamics} observe that the rate of new infections in node $i$ is described by the quantity $f_i = a_ix_i \bar y$. Moreover,
		$$
		\frac{\dot f_i}{f_i} = \frac{\dot{\bar y}}{\bar y} + \frac{\dot x_i}{x_i} = \tilde x - \gamma - a_i\bar y = w_i.
		$$
		Hence, $w_i$ is the logarithmic derivative of the rate of new infections in node $i$.}
Let also
\be\label{hattdef}\hat t=\inf\{t\ge0:\, \tilde x(t)\le\gamma\}\,\ee
\textcolor{black}{denote the peak time of the weighted aggregate of infected individuals $\bar{y}$,} and observe that Theorem \ref{theo:bary} implies that $\hat t<+\infty$. 
Also, for every $i=1,\ldots, n$, let 
$$\ov t_i=\inf\{t\!\ge\!0:\, w_i(t)\!\le\!0\}=\inf\{t\!\ge\!0:\, \tilde x(t)\!\le\!\gamma+a_i\ov y(t)\}$$
\textcolor{black}{be the first time instant at which the rate of new infections in node $i$ starts decreasing,}
and notice that \be \label{ov_t}\ov t_i\le\hat t\,,\ee
and $\ov t_i\le \ov t_j$ if and only if $a_j\le a_i$.  Hence, it is possible to order these time instants from the entries of vector $a$.

The following are two useful technical results.  

\begin{lemma}\label{lemma2new}
Consider the rank-$1$ network SIR epidemic model \eqref{eq:dynamics} and 
 every initial state $(x(0),y(0))$ such that $y(0)\gneq\0$. 
Then, for every $i=1,\ldots, n$, 
\begin{enumerate}
\item[(i)] 
$\dot w_i(t)<-a_i\ov yw_i(t)$ for every $t\ge0$; 
\item[(ii)] 
$w_i(t)$ is strictly decreasing for $0\le t\le\ov t_i$; 
\item[(iii)] $w_i(t)<0$ for every $t> \ov t_i$;
\item[(iv)] for every $t\ge0$,  
\be\label{ddotyi}\ddot y_i(t)=a_ix_i(t)\ov y(t)w_i(t)-\gamma\dot y_i(t)\,;\ee  
\item[(v)] if $\dot y_i(t)=0$ for some $t\ge\ov t_i$, then $t$ cannot be a local minimum point of $y_i(t)$.  
\end{enumerate}
\end{lemma}\smallskip
\begin{proof} 
See Appendix \ref{sec:proof-lemma2new}.
\end{proof}\medskip

\begin{proposition}\label{prop:main}
Consider the connected rank-$1$ network SIR epidemic model \eqref{eq:dynamics} \textcolor{black}{with initial state $(x(0),y(0))$ in $\mc S$ such that $y(0)\gneq\0$.}
   Then, for every $i=1,\ldots,n$, 
\begin{enumerate}
\item[(i)] $y_i(t)$ admits at most one local minimum time $\check t_i\ge0$.
\end{enumerate}
Moreover, if such local minimum time $\check t_i$ exists, 
\begin{enumerate}
\item[(ii)]  it satisfies
\be \label{ti}0\le \check t_i\le\ov t_i\,,\ee
with $\check t_i=\ov t_i=0$ if and only if $w_i(0)\le0$ and $\dot y_i(0)>0$; 
\item[(iii)] it cannot occur after any stationary local maximum point of $y_i(t)$.
\end{enumerate}
\end{proposition}\smallskip
\begin{proof} 
If $w_i(0)\le0$, then $\ov t_i=0$ and Lemma \ref{lemma2new}(v) implies that no stationary point $t\ge0$ of $y_i(t)$ can be a local minimum point. It follows that the only local minimum point of $y_i(t)$ can possibly be $\check t_i=0$, which is the case if and only if $\dot y_i(0)>0$. Hence the claims are proved in the special case $w_i(0)\le0$. 

On the other hand, if $w_i(0)>0$, then $\ov t_i>0$ and the interior extremum theorem and Lemma \ref{lemma2new}(v) imply that there cannot be any minimum points of $y_i(t)$ in the interval $[\ov t_i,+\infty)$. 
This proves point (ii). We are then left with studying local minimum points of  $y_i(t)$ in the interval $[0,\ov t_i)$. 
%
Let $s\ge0$ be a stationary local maximum point of $y_i(t)$, and let $u$ in $(s,\ov t_i)$ be a (necessarily stationary) local minimum point of $y_i(t)$. 
Then, \be\label{dotyius=0}\dot y_i(s)=\dot y_i(u)=0\,,\ee and \be\label{ddotyius>=0}
\ddot y_i(s)\le0\,,\qquad \ddot y_i(u)\ge0\,.\ee 

\tcb{Notice that we cannot have $y_i(s)=0$ or otherwise $y_i(t)=0$ in a neighborhood of $s$ (as $s$ is a local maximum point for $y_i(t)$) thus contradicting Proposition \ref{lemma:basic}(v).}  Hence, we get
\be\label{ddotys}\ba{rcl}
0&\ge&\ddot y_i(s)/y_i(s)\\[5pt]
&=&a_ix_i(s)\ov y(s)w_i(s)/y_i(s)-\gamma\dot y_i(s)/y_i(s)\\[5pt]
&=&\gamma w_i(s)\\[5pt]
&>&\gamma w_i(u)\\[5pt]
&=&a_ix_i(u)\ov y(u)w_i(u)/y_i(u)-\gamma\dot y_i(u)/y_i(u)\\[5pt]
&=&\ddot y_i(u)/y_i(u)\\[5pt]
&\ge&0
\,,\ea\ee
where the first and the last inequalities above follow from \eqref{ddotyius>=0}, the first and the last 
identities follow from \eqref{ddotyi}, 
the other two identities from \eqref{dotyius=0} and the fact that, by the second equation in \eqref{eq:dynamics}, $a_{i} x_{i}(t)\bar{y}(t)=\gamma y_{i}(t)$ when $\dot y_i(t)=0$,  
and the strict inequality in the  middle holds true because of Lemma \ref{lemma2new}(ii). 
As \eqref{ddotys} is a contradiction, this shows that a local minimum point \tcb{$u<\ov t_i$} of $y_i(t)$ cannot follow any stationary local maximum point $s\ge0$ of $y_i(t)$, thus proving point (iii). 

Finally, to prove point (i), assume by contradiction that there exist two distinct local minimum points $r<u$ of $y_i(t)$ in the interval \tcb{$[0,\ov t_i)$}. 
Then, there would necessarily exist a local maximum point of $y_i(t)$ in the interval $(r,u)$. But, since $s>r\ge0$, such local maximum point would also be stationary, thus violating point (iii). 
Therefore, there cannot exist two distinct local minimum points $r<u$ of $y_i(t)$ in the interval \tcb{$[0,\ov t_i)$}, thus completing the proof of point (i). 
\end{proof}\medskip

\begin{remark}\label{remark:peak_time}
\tcb{Equations \eqref{ti} and \eqref{ov_t} imply that the local minimum point of $y_i$ can never occur after the peak of the weighted aggregate of infected individuals $\bar{y}$.    
Hence, if at some time $\tau\ge0$ both $\dot{\ov y}(\tau)\le0$ (so that $\ov y$ has peaked at some time $\hat t\le\tau$) and $\dot y_i(\tau)<0$ (so that the epidemic in $i$ is currently regressing), then $y_i(t)$ will remain decreasing for all $t\ge\tau$.}
\end{remark}\medskip

As a consequence of Proposition \ref{prop:main}, we get the following result classifying the possible dynamic behaviors of the fraction of infected individuals in the single populations of the network SIR epidemic model with rank-$1$ interaction matrix. This classification is based on the study of the sign of two quantities:
\begin{equation*} 
	\dot y_i(0) = a_i x_i(0)\ov y(0)-\gamma y_i(0), \, w_i(0)=\tilde x(0)-\gamma-a_i\ov y(0) .
\end{equation*}

\begin{theorem}\label{coro:main}
Consider the connected rank-$1$ network SIR epidemic model \eqref{eq:dynamics} \textcolor{black}{with initial state $(x(0),y(0))$ in $\mc S$ such that $y(0)\gneq\0$.}
Then,  for every $i=1,\ldots,n$, 
\begin{enumerate}
\item[(i)] if 
\be\label{ai<=gamma}\dot y_i(0)\le0\,,\ee 
and 
\be\label{subsubcritical}w_i(0) \le 0\,,\ee 
then $y_i(t)$ is strictly decreasing for $t\ge0$; 
\item[(ii)] if \be\label{ai>gamma}\dot y_i(0) >0\,,\ee 
or if \be\label{ai=gamma}\dot y_i(0) =0\,,\ee 
and \be\label{subsupercritical}w_i(0)>0\,,\ee 
then there exists a peak time $\hat t_i>0$ such that $y_i(t)$ is strictly increasing on $[0,\hat t_i]$ and strictly decreasing on $[\hat t_i,+\infty)$;
\item[(iii)] if 
\be\label{ai<gamma}\dot y_i(0)<0\,,\ee   and 
\be\label{subsupercritical2}w_i(0)>0\,,\ee
then either $y_i(t)$ is strictly decreasing for $t\ge0$ or there exist a local minimum time $\check t_i$ and a peak time $\hat t_i$ such that $0<\check t_i<\hat t_i$ and $y_i(t)$ is strictly decreasing on $[0,\check t_i]$, strictly increasing on $[\check t_i, \hat t_i]$, and strictly decreasing on $[\hat t_i,+\infty)$;
\end{enumerate}
\end{theorem}
\smallskip
\begin{proof} 
(i) If \eqref{subsubcritical} holds true, then $\ov t_i=0$. On the other hand, \eqref{ai<=gamma} and Proposition \ref{prop:main}(ii) rule out the possibility that there exists any minimum point for $y_i(t)$. Therefore, $y_i(t)$ is strictly decreasing for $t\ge0$. 
 
(ii) If \eqref{ai>gamma} holds true, then Proposition \ref{prop:main}(i) implies that $\check t_i=0$ is the only minimum point of $y_i(t)$. 
On the other hand, if \eqref{ai=gamma} and \eqref{subsupercritical} both hold true, then it follows from \eqref{ddotyi} that 
$$\ddot y_i(0)=a_ix_i(0)\ov y(0)w_i(0)-\gamma\dot y_i(0)=\gamma\ov y^2(0)w_i(0)>0\,,$$
(where the second identity follows from the fact that, by the second equation in \eqref{eq:dynamics}, $a_{i} x_{i}(t)\bar{y}(t)=\gamma y_{i}(t)$ when $\dot y_i(t)=0$)
thus implying that also in this case $\check t_i=0$ is a local minimum point for $y_i(t)$. 

\tcb{Since $y(0)\gneq\0$}
and, by Proposition \ref{prop:main}(i), $y_i(t)$ cannot have another local minimum points besides $\check t_i=0$, 
it follows that  exists a peak time $\hat t_i>0$ such that $y_i(t)$ is strictly increasing for $t$ in $[0,\hat t_i]$ and strictly decreasing for $t$ in $[\hat t_i,+\infty)$. 

(iii) From  \eqref{ai<gamma}, $0$ is a nonstationary local maximum point of $y_i(t)$. Since, by Proposition \ref{prop:main}(i), $y_i(t)$ can have at most one local minimum point, and 
by Proposition \ref{lemma:basic}(vi), 
\be\label{limyi=0}\lim_{t\to+\infty}y_i(t)=0<y_i(0)\,,\ee
it follows that either $y_i(t)$ is strictly decreasing for $t\ge0$ (in case there is no local minimum point) or, if a local minimum point $\check t_i>0$ exists, then there exists also a peak time $\hat t_i>\check t_i$ so that $y_i(t)$ is strictly increasing on $[0,\hat t_i]$ and strictly decreasing on $[\hat t_i,+\infty)$.
\end{proof}\medskip

%

The previous result provides a classification of the behavior of \tcb{infection curves of} the single subpopulations. 
In particular, observe that if $\dot y_i(0)=0$, then the behavior of the single \tcb{infection curves} depends only on the sign of $w_i(0)$: in this case, Theorem \ref{coro:main} provides a tight condition.

In Theorem \ref{coro:main}, each condition is considered from the perspective of the single subpopulation $i$. However, note that if \eqref{ai<=gamma} is true for all $i=1,\ldots,n$, meaning that the \tcb{infection curves are initially decreasing, then they will decrease forever}. Indeed, multiplying both sides of \eqref{ai<=gamma} by $b_i$  and summing over $i$, we would get $ \tilde x(0)<\gamma$ so that $w_i(0)<0$ for all $i$.

\subsection{Sufficient conditions for multimodal infection curves}\label{sec:sufficient}
We now consider a particular class of rank-$1$ interaction matrices in the form
\be\label{a_spec} A = \beta \1 b^T\,, \ee
with $\beta >0$ and $\1^T b=1$. This is a special case of the one studied in Section \ref{sec:4} where the vector $a$ has all equal entries and the entries of $b$ sum up to $1$. 

\begin{remark}
	 This model corresponds to a scenario in which all individuals have the same susceptibility to the disease but different \tcb{infectivity level}. E.g., individuals wearing medical masks become infected with the same probability but spread the disease differently. Note that a simple case of this class of matrices is $A=\1\1^T$, studied in Example \ref{example:mutimodal}.	The network SIR epidemic model with this interaction matrix is of interest for control applications, c.f.~\cite{acemoglu2021optimal}. Indeed, even if the dynamics at the nodes are homogeneous and thus the infection spreads at the same rate, it may be convenient to divide individuals into multiple groups, e.g., to study the effects of differentiated control policies, especially in cases whereby the cost of applying a control and epidemic cost for the diffusion of the disease may differ depending on the age of the individuals. 
\end{remark}

We observe that, for rank-$1$ interaction matrices in the form \eqref{a_spec}, the dynamics become
\be\label{SIR-single-special}\dot{x}_{i} = -\beta x_{i} \bar{y}, \qquad \dot{y}_{i} = \beta x_{i} \bar{y} - \gamma y_i\,,\ee 
for every $i=1,\dots,n$, and
\be\label{SIR-specialcase}\dot{\ov x}=-\beta \ov x\,\ov y, \qquad \dot{\ov y}=\ov y\left(\beta \ov x-\gamma\right)\,,\ee
since $\ov x$ and $\tilde{x}$ differ in a constant term only. The next result provides sufficient conditions for multimodality of the infection curve at the single node level and encompasses Example \ref{example:mutimodal}.
We first need to define auxiliary functions
\be \label{g(eps)} g_i(\eps) = \frac{1-\eps}{1-b_i\eps}\left(1-\frac{\gamma}{\beta} + \frac{\gamma}{\beta} \log\frac{\gamma}{\beta(1-b_i\eps)} \right)-\eps\,.
	\ee
Notice that $$g_i(0)=1-\frac{\gamma}{\beta} + \frac{\gamma}{\beta} \log\frac{\gamma}{\beta},\qquad g_i(1)=-1\,.$$
As a consequence, when $\gamma/\beta <1$, $g_i$ admits zeroes in $[0,1]$ and we put
\be \label{oveps2}\ov \eps_i=\min\left\lbrace \eps\in[0,1]:\, g_i(\eps)=0\right\rbrace.\ee
%
%
\begin{proposition}\label{prop:special}
	Consider the rank-$1$ network SIR epidemic model \eqref{eq:dynamics} with $a =\beta \1$ and $\1^T b=1$.
	Consider a subpopulation $i$ in $\{1,\ldots,n\}$ and an initial state $(x(0), y(0))$ in $\mc S$ that satisfy the following conditions:
	\begin{eqnarray}
	\label{no_recovered}  x(0)+y(0) = \mb 1\,,\\
	\label{1st_assump} \beta x_i(0)\bar{y}(0)-\gamma y_i(0) <0\,, \\
	\label{2nd_assump} \beta \ov x(0)>\gamma\,, \\
	\label{3nd_assump} 0<y_i(0)<\ov \eps_i\,.
	\end{eqnarray}
	Then, there exist a local minimum time $\check t_i$ and a peak time $\hat t_i$ such that $0<\check t_i<\hat t_i$ and $y_i(t)$ is strictly decreasing on $[0,\check t_i]$, strictly increasing on $[\check t_i, \hat t_i]$, and strictly decreasing on $[\hat t_i,+\infty)$.
\end{proposition}\smallskip
\begin{proof}
	From \eqref{1st_assump} and \eqref{SIR-single-special} it follows that $\dot y_i(0)<0$, which implies that $y_i(t)$ is strictly decreasing for sufficiently small $t>0$. On the other hand, \eqref{SIR-specialcase} and \eqref{2nd_assump} imply that $\bar y(t)$ is strictly increasing for sufficiently small $t>0$. Since $\ov x$ and $\ov y$ satisfy the scalar autonomous SIR epidemic model \eqref{SIR-specialcase}, this implies that $\ov y(t)$ has a peak at some time $\hat t>0$ and 
	\be \label{eq:peak} \ov x (\hat t)= {\gamma}/{\beta}\,.\ee 
	From Proposition \ref{prop:scalar-SIR}(i) we obtain that the peak value of the weighted aggregate of infected individuals is
	\be \label{y_max} \ba{rcl}
		\ov y(\hat t\,) &=&\ds \ov x(0) + \ov y(0) -\ov x(\hat t) + \frac{\gamma}{\beta} \log \frac{\ov x(\hat t)}{\ov x(0)}\\[2pt]
		&=&\ds 1 -\frac{\gamma}{\beta} + \frac{\gamma}{\beta} \log \frac{\gamma}{\beta \ov x(0)} \,,
	\ea\ee
	where the second equality follows from \eqref{no_recovered} and \eqref{eq:peak}.
	Moreover, \eqref{SIR-single-special} and \eqref{SIR-specialcase} imply that
	${x_i(t)}/{x_i(0)} = {\ov x(t)}/{\ov x(0)}$, 
	for  $i=1,\dots,n,$ and $t\ge0$. Therefore, using \eqref{eq:peak} we obtain
	\be\label{xi_ovx_hatt} x_i(\hat t) = \frac{\gamma}{\beta}\frac{x_i(0)}{\ov x(0)} \,,\ee
	We now prove that $y_i(t)$ cannot remain decreasing for all $t>0$. 
	Assume by contradiction that 
	$\dot y_i(t)\le 0$ for all $t$ in $[0,\hat t]$. 
	In particular, this implies that $y_i(\hat t)\le y_i(0)$.
	This together with \eqref{xi_ovx_hatt} and \eqref{y_max} imply that
	\be\label{ineq1} \ba{rclcl}
	0 &\geq&  \dot y_i(\hat t\,)\\[2pt]
	& =& \ds \beta x_i(\hat t\,) \Bar{y}(\hat t\,) - \gamma y_i(\hat t\,) \\[2pt]
	& =& \ds \gamma \frac{x_i(0)}{\ov x(0)}\left( 1-\frac{\gamma}{\beta}+ \frac{\gamma}{\beta}\log \frac{\gamma}{\beta \ov x(0)}\right)- \gamma y_i(\hat t\,) \\[5pt]
	& \geq& \ds \gamma \left[ \frac{x_i(0)}{\ov x(0)}\left( 1-\frac{\gamma}{\beta}+ \frac{\gamma}{\beta}\log\frac{\gamma}{\beta \ov x(0)}\right)-y_i(0) \right]\,.
	\ea \ee
	Notice that, because of the assumptions on $b$, we have that $\ov x(0)= 1-\ov y(0) \leq 1- b_i y_i(0)$. Since the last expression in \eqref{ineq1} is decreasing in $\ov x(0)$ and $x_i(0)=1-y_i(0)$, we get that 
	\be \label{eq:g_i}  
	 \ds \gamma g_i (y_i(0))\le 0\,. 
	\ee
	By \eqref{2nd_assump}, we necessarily have that $\gamma/\beta<1$. This implies that $g_i(0)>0$ and, with \eqref{eq:g_i}, that $y_i(0)\ge\ov\eps_i$, thus violating \eqref{oveps2}. This contradiction implies that  $y_i(t)$ cannot remain decreasing for all $t>0$. The claim then follows from Theorem \ref{coro:main}.
\end{proof}\medskip

\begin{remark}
	Observe that the set of model parameters and initial states that satisfy the assumptions of Proposition \ref{prop:special} is nonempty. To prove this, consider a network with $n$ nodes, interaction matrix as in \eqref{a_spec} with parameters $\beta >\gamma$ and
	\be \label{eq:b1} b_1 < \min\left\{\frac{\gamma}{\beta}, 1-\frac{\gamma}{\beta}\right\}\,, 
	\ee 
	Fix an initial state $(x(0),y(0))$ in $\mc S$ such that $0<y_1(0)<\ov\epsilon_1$ and $y_j(0)=0$ for $2\le j\le n$. Notice that $\ov y(0) = b_1 y_1(0)$. A straightforward check shows that \eqref{1st_assump} and \eqref{2nd_assump} are automatically satisfied putting no further restriction on $y_1(0)$. Therefore, all assumptions of Proposition \ref{prop:special} are satisfied.
%
%
%
\end{remark}\medskip

\begin{remark}
	Observe that under the assumptions of Proposition \ref{prop:special} we can provide an upper bound for stationary infection peaks of a node $i$. Let $\hat t_i$ be the peak time of node $i$, so that $\dot{y}_i(\hat t_i)=0$. By \eqref{SIR-single-special}, this implies that 
	\be y_i(\hat t_i) = \frac{\beta}{\gamma} x_i(\hat t_i) \ov y(\hat t_i) \,.\ee
	Since the weighted aggregate of infected individuals is limited above by its infection peak value, i.e. $\ov y(t)\leq\ov y(\hat t\,)$ for $t \geq 0$, and the fraction of susceptible individuals is monotonically decreasing, we have that 
	$$ \ba{rclcl}
	y_i(\hat t_i) &\leq & \ds \frac{\beta x_i(0)}{\gamma}  \ov y(\hat t) \\[2pt]
	&=& \ds \frac{\beta x_i(0)}{\gamma}\left(\ov x(0) + \ov y(0) -\ov x(\hat t) + \gamma \log \frac{\ov x(\hat t\,)}{\ov x(0)}\right) \\[2pt]
	&=& \ds x_i(0)\left(\frac{\beta}{\gamma} -1 +\log \frac{\gamma}{\ov x(0)}\right)\,,
	\ea $$
	where the first equivalence follows from \eqref{y_max} and the last one from \eqref{eq:peak} and \eqref{no_recovered}.
	
\end{remark}\medskip

\section{Numerical Simulations}\label{sec:5}
In this section, we present the result of some numerical simulations of the network SIR epidemic model. Some more numerical simulations are reported in \cite{9992408}. 

\begin{figure}
	\centering
	\includegraphics[scale=0.7]{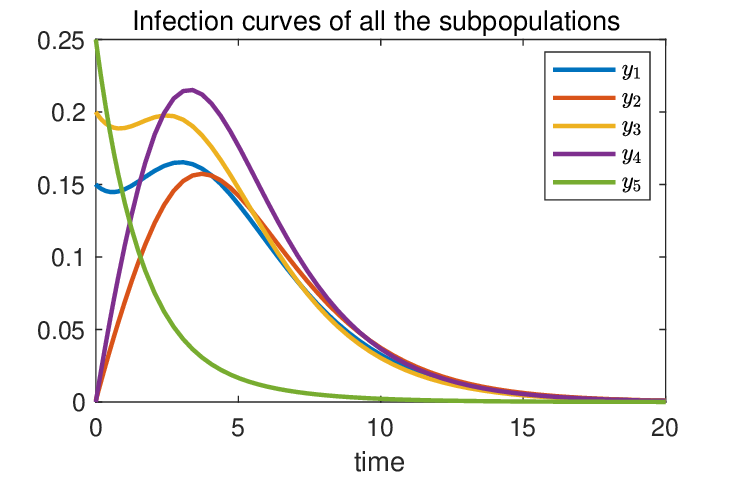}\\[13pt] 
	\includegraphics[scale=0.7]{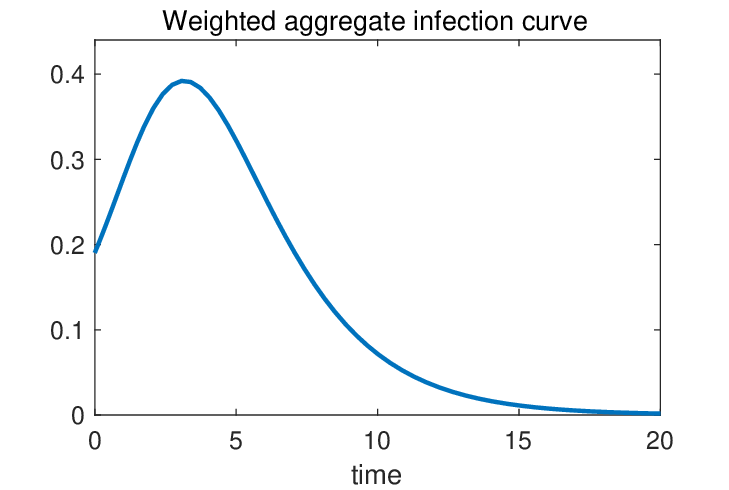}
	\caption[]{Numerical simulations of the network SIR epidemic model with $n=5$ nodes and rank-$1$ interaction matrix. }	
	\label{fig:fig2}
\end{figure}

We start by reporting the infections curves of a network SIR epidemic model with $n=5$ nodes and rank-$1$ interaction matrix. \textcolor{black}{In Figure \ref{fig:fig2}, the interaction matrix is $A = ab^{T}$ with $a = (0.1, 0.25, 0.6, 1, 0.2)$ and $b = (0.45, 0.4, 0.6, 0.65, 0.01)$, and the recovery rate is $\gamma= 0.6$. The initial state has entries $x_{1}(0)= 0.85$, $x_{2}(0)= 0.999$, $x_{3}(0)= 0.8$, $x_{4}(0)= 1$, $x_{5}(0)= 0.75$, and $y(0)=\1-x(0)$. }The top plot in Figure \ref{fig:fig2} shows that the infection curves of nodes $1$ and $3$ are both multimodal with two changes of monotonicity: for both $i=1$ and $i=3$, the fraction of infected individuals $y_i(t)$ is strictly decreasing for times $t$ in $[0,\check t_i]$, it is strictly increasing in $[\check t_i,\hat t_i]$, and it is strictly decreasing for $t$ in $[\hat t_i,+\infty)$. The bottom plot displays the weighted aggregate infection curve $t\mapsto\ov y(t)$, which is unimodal  with a peak in $\hat t$, as proved in Theorem \ref{theo:bary}. Notice also that, as consequence of Proposition \ref{prop:main}(ii) (cf.~Remark \ref{remark:peak_time}), $\check t_i \le \hat t$ for  $i = 1,3$, namely the local minimum of the infection curve in each of the two nodes cannot occur after the  peak time of the weighted aggregate of infected individuals. 

\begin{figure}
	\centering
	\includegraphics[scale=0.7]{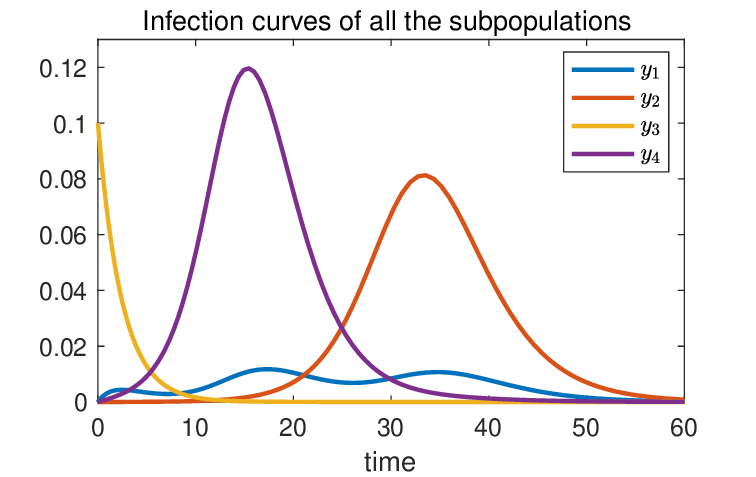}\\[10pt]
	\includegraphics[scale=0.7]{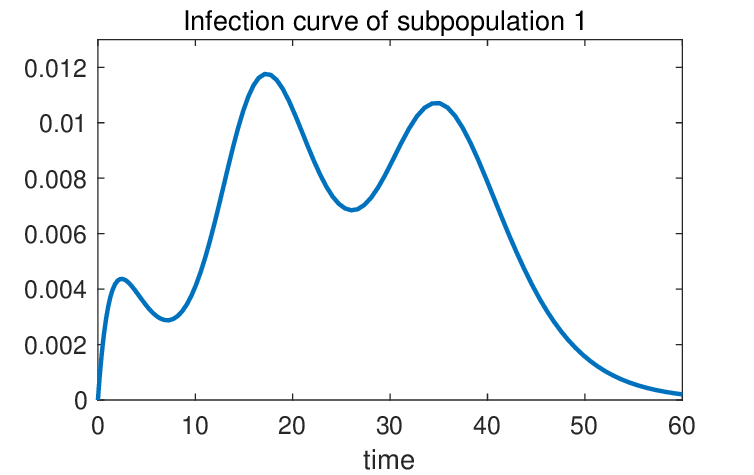}
	\caption[]{\label{fig:fig5}Numerical simulation of the network SIR epidemic model with $n=4$ nodes and full-rank interaction matrix. }	
\end{figure}
\textcolor{black}{In Figure \ref{fig:fig5}, we report numerical simulations of a network SIR epidemic model with $n=4$ nodes, recovery rate  $\gamma= 0.5$, and interaction matrix 
	\begin{equation*}
		A=\left(\ba{cccc}0.05& 0.07& 0.05& 0.05\\ 0.0001& 0.8& 0.0001& 0.0001\\  0.0001&	0.0001& 0.1& 0.0001\\ 0.01& 0.01& 0.01& 0.9\ea\right)\,.
	\end{equation*}
The initial state has entries $x_1(0)=x_2(0)=x_4(0)=1$, $y_1(0)=y_2(0)=y_4(0)=0$, and $x_{3}(0)=0.9=1- y_{3}(0)$, i.e., all subpopulations are initially completely susceptible except for subpopulation $3$ that has $10\%$ of infected individuals. 
It can be observed that the infection curve in the first node displays three peaks. Moreover, we observe the presence of delays among the peaks in different nodes. This simulation shows that the limitation of the number of peaks to two is a peculiar feature of rank-$1$ interaction matrices, while for general networks, even with a limited number of nodes, multiple peaks may emerge. This is an interesting feature of the network SIR epidemic model especially since actual epidemic data often displays multimodal infection curves with multiple peaks occurring at different times in different subpopulations. (See, e.g.,  \cite{covid-cases} for COVID-19 infection curves in the different countries.)}


\section{Conclusion}
\label{sec:conclusion}
In this paper, we  have studied the network SIR epidemic model, focusing on the special case of rank-$1$ interaction matrices. \tcb{We have identified $n$ invariants of motion and built on them to find explicit expression for the limit equilibrium point as a function of the initial state and to analyze the stability of the equilibrium points}. We have then proved that, in contrast to the scalar SIR model, in the network SIR epidemic model the infection curve associated to a single subpopulation may be multimodal, and we have established sufficient conditions for the occurrence of this phenomenon. Furthermore, we have characterized all the possible behaviors that the dynamics may exhibit at single node level, showing that the infection curve in a single node can undergo two changes of monotonicity at most.  Finally, we have conducted a numerical analysis showing that for more general interaction matrices the network SIR epidemic model may exhibit more than two peaks at single node level. 

 \tcb{We are aware that the phenomenon of multiple waves of infection cannot be fully explained by the heterogeneity introduced by the network and is also largely determined by the adaptive behavior and endogenous response of individuals to the epidemic, as well as by the phenomenon of waving immunity \cite{Ehrhardt2019}. For example, some papers have studied models that take into account how individuals adapt their behavior, resulting in a modification of the parameters of the model at macroscopic level \cite{CAPASSO197843, Alutto2021OnSE}, or in relation to loss of immunity over time \cite{Jin2007, Li2014}. Other works study the epidemics dynamics coupled with an evolutionary game-theoretic decision-making mechanism concerning some individual behaviors \cite{frieswijk2022meanfield}, \cite{elokda2021dynamic}.
We also acknowledge the many possible extensions of the network SIR epidemic model to more than three compartments to keep track for instance of the many forms of infection and possibly vaccination
\cite{Li1995, Giordano.ea:2020, Birge.ea:2020}. 
Future work aims to include in the model more complex phenomena to fully describe the occurrence of multiple peaks in epidemic models.}

\bibliographystyle{IEEEtran}
\bibliography{bib}

\appendices
\section{Proof of Lemma \ref{lemma1}}\label{sec:proof-lemma1}
By taking the time derivative of both sides of the first equation in \eqref{barxy} and substituting the first one in \eqref{eq:dynamics}, we get 
$$\dot{\ov x}=\sum_{j=1}^{n} b_{j}\dot x_{j}=-\bar{y}\sum_{j=1}^{n}b_ja_{j} x_{j}=-\ov y\tilde x\,,$$
thus proving the first equation in \eqref{dotovx+dotovy}. Analogously, by taking the time derivative of both sides of the second equation in \eqref{barxy} and substituting the second equation in \eqref{eq:dynamics}, we get 
$$\dot{\ov y}=\sum_{j=1}^{n} b_{j}\dot y_{j}= \bar{y}\sum_{j=1}^{n}b_ja_{j} x_{j} - \gamma \ov y=\ov y\left(\tilde x-\gamma\right)\,,$$
thus proving the second equation in \eqref{dotovx+dotovy}.

\section{Proof of Lemma \ref{lemma:varphi}}\label{sec:proof-lemmavarphi}
\tcb{For $x=\0$, equation \eqref{eq:barx*} reduces to $\xi=0$, which has a unique solution $\varphi(\0,y)=0=\ov x$. 
This proves claims (i)--(ii) 
in the special case $x=\0$ and the ``if'' part of claim (iii). 
For $x\gneq\0$, consider the function $g:\R_+\to\R$ defined by
$$ g(\xi) =   \sum_{j=1}^n b_j c_j \exp (a_j \xi / \gamma)-\xi\,,$$ 
where $c_j=x_j\exp(-a_j(\ov x+\ov y)/\gamma)$,  for $1\le j\le n$.  Observe that equation \eqref{eq:barx*} is equivalent to $g(\xi)=0$. 	
We have that:
\begin{enumerate}
\item[(a)] $g(\xi)$ is differentiable and strictly convex in $\xi\ge0$;   
\item[(b)] $g(0)=\sum_jb_jc_j>0$ and $g(\xi)\stackrel{\xi\to+\infty}{\longrightarrow} +\infty$;  
\item[(c)] $g(\ov x)=\sum_jb_jx_j\exp(-a_j\ov y/\gamma)-\ov x\le\sum_jb_jx_j-\ov x =0$,
with equality if and only if $y=\0$; 
\item[(d)] if $y=\0$ then $g'(\ov x)=\sum_ja_jb_jx_j/\gamma-1=\tilde x/\gamma-1<0$. 
\end{enumerate}
Points (a)--(c) imply that $g(\xi)$ has at least one and at most two zeros. 
By (b), $\varphi(x,y)>0$ for $x\gneq\0$, proving the ``only if'' part of (iii). 
For $y\gneq\0$, by points (b)--(c), $g(\xi)$ has one zero in $(0,\ov x)$ and another one in $(\ov x,+\infty)$, proving claim (i) and the ``only if'' part of claim (ii). 
For $y=\0$, (b)--(d) imply that $g(\xi)$ has no zeros in $[0,\ov x)$ and $g(\ov x)=0$, proving claims (i) and the ``if'' part of claim (ii). 
Finally, since $g'(\varphi(x,y))<0$ for $y\gneq0$, the Implicit Function Theorem implies that $\varphi(x,y)$ is differentiable for $(x,y)$ in $\mc S$ such that $x\gneq\0$. Continuity of $\varphi(x,y)$ on the whole space $\mc S\setminus\{(x,\0):\tilde x\ge\gamma\}$ then follows upon noting that 
$0\le \varphi(x,y)\le\ov x\stackrel{x\to\0}{\longrightarrow}0$, so that $\varphi(x,y)\stackrel{x\to\0}{\longrightarrow}0=\varphi(\0,y)$. }

\section{Proof of Lemma \ref{lemma2new}}
\label{sec:proof-lemma2new}
(i) From the first equation \eqref{eq:dynamics}, we have that 
\be\label{dottildex}\dot{\tilde x}=\sum_ia_ib_i\dot x_i=-\ov y\sum_ia_i^2b_ix_i\le0\,.\ee
By Proposition \ref{lemma:basic}(v), the assumption $y(0)\gneq\0$ implies that 
\be\label{y(t)>0}y(t)>0\,,\qquad\forall t>0\,,\ee
so that in particular 
\be\label{ovy>0}\ov y(t)=\sum_jb_jy_j(t)>0\,,\qquad \forall t\ge0\,.\ee 
 It then follows from \eqref{widef},\eqref{dotovx+dotovy}, \eqref{dottildex}, and \eqref{ovy>0} that 
$$\dot w_i
=\dot{\tilde x}-a_i\dot{\ov y}
=\dot{\tilde x}-a_i\ov y(\tilde x-\gamma)
=\dot{\tilde x}-a_i\ov yw_i-a_i^2\ov y^2
<-a_i\ov yw_i\,,$$
thus proving the claim. 

(ii) For $0\le t\le\ov t_i$, we have $w_i(t)\ge0$ so that, by point (i) and \eqref{ovy>0}, we have that 
$\dot w_i(t)<-a_i\ov y(t)w_i(t)\leq0\,.$
This implies that $w_i(t)$ is strictly decreasing for $0\le t\le\ov t_i$. 

(iii) It follows from (i) that, for every $t_i^*\ge0$ such that $w_i(t_i^*)=0$, we have 
$\dot w_i(t_i^*)<-a_i\ov y(t_i^*)w_i(t_i^*)=0\,,$
so that $w_i(t)<0$ for $t>t_i^*$. 
If $w_i(0)\le0$ so that $\ov t_i=0$, this implies that  $w_i(t)<0$ for all $t>0=\ov t_i$. 
On the other hand, $w_i(0)>0$ so that $\ov t_i>0$,  it follows that $w_i(t)<0$ for all $t>\ov t_i$. 
%

(iv) Taking the derivative of both sides of the second line in \eqref{eq:dynamics}  and substituting  the first line of \eqref{eq:dynamics} and \eqref{dotovx+dotovy} yield
$$
\ba{rclcl}\ddot y_i&=& a_{i}\left(\dot x_{i} \bar{y}+x_{i} \dot{\bar{y}}\right) - \gamma\dot y_{i}\\[3pt]
&=&a_ix_i\ov y\left(\tilde x-\gamma-a_i\ov y\right) - \gamma\dot y_{i}
&=&a_ix_i\ov y w_i-\gamma\dot y_i\,,\ea$$
thus proving \eqref{ddotyi}.

(v) 
Assume by contradiction that $t\ge\ov t_i$ is a  local minimum point of $y_i(t)$ with $\dot y_i(t)=0$. By \eqref{ddotyi}, we then  have 
\be\label{ddotyit}\ddot{y}_i(t) 
= a_ix_i(t)\ov y(t) w_i(t)-\gamma\dot y_i(t)\\[5pt]
=\gamma y_i(t)w_i(t)\,,\ee
where the last identity holds true since, by the second equation in \eqref{eq:dynamics}, $\dot y_i(t)=0$  is equivalent to  \be\label{equivalent}a_{i} x_{i}(t) \bar{y}(t)=\gamma y_{i}(t)\,.\ee
Equations \eqref{ddotyit} and \eqref{y(t)>0} imply that 
\be\label{sgnddotyi}\sgn(\ddot{y}_i(t))=\sgn(w_i(t))\,.\ee
Now, recall that by point (iii) we have $w_i(t)\le0$ for every $t\ge\ov t_i$. If $w_i(t)<0$, then equation \eqref{sgnddotyi} implies that $\ddot y_i(t)<0$, so that $t$ cannot be a local minimum point for $y_i(t)$. 
On the other hand, if $w_i(t)=0$, by point (i) and \eqref{ovy>0} we get that 
\be\label{dotwi<0}\dot w_i(t)<-a_i\ov y(t)w_i(t)=0\,,\ee
while \eqref{sgnddotyi} implies that $\ddot y_i(t)=0$. Hence, taking the derivative of both sides of \eqref{ddotyi}, by \eqref{equivalent} and $w_i(t)=0$, we get 
$$\ba{rcl}\dddot y_i(t)&=&a_ix_i(t)\ov y(t)\dot w_i(t)+a_i(\dot{(x_i(t)\ov y(t))})w_i(t)-\gamma\ddot y_i(t)\\[5pt]
&=&\gamma y_{i}(t)\dot w_i(t)\\[5pt]
&<&0\,,
\ea$$
where the last inequality follows from \eqref{y(t)>0} and \eqref{dotwi<0}. Together with $\dot y_i(t)=\ddot y_i(t)=0$, the above implies that $t$ is an inflection point for $y_i(t)$, hence it particular it is not a local minimum point. We have thus proven that $y_i(t)$ cannot have any stationary minimum points in the interval $[\ov t_i,+\infty)$.

\newpage
\begin{IEEEbiography}[{\includegraphics[width=1in,height=1.25in,clip,keepaspectratio]{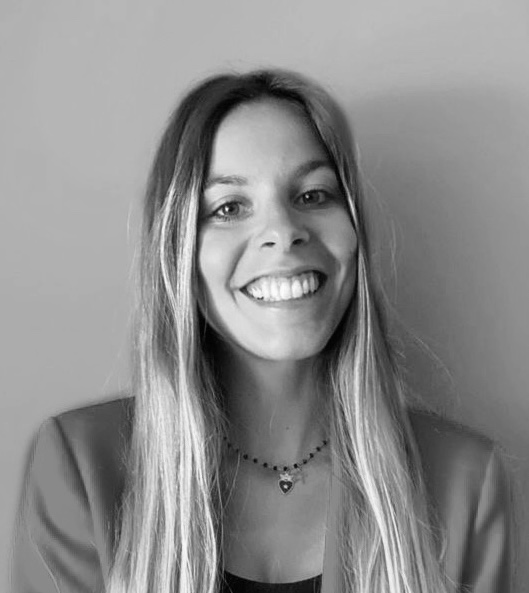}}]{Martina Alutto} received the B.Sc.~and the M.S.~(cum laude) in Mathematical Engineering from  Politecnico  di  Torino, Italy, in  2018 and 2021, respectively. She is currently a PhD student in Pure and Applied Mathematics at the Department of Mathematical Sciences, Politecnico di Torino, Italy. Her research interests focuses on analysis and control of network systems, with application to epidemics and social networks.
\end{IEEEbiography}

\begin{IEEEbiography}[{\includegraphics[width=1in,height=1.25in,clip,keepaspectratio]{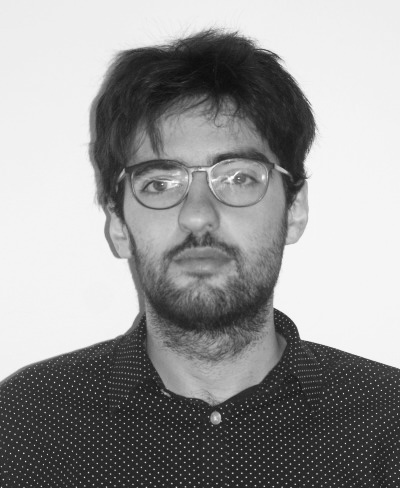}}]{Leonardo Cianfanelli} received the B.Sc.  in Physics and Astrophysics in 2014 from Universit\`a  di Firenze, Italy, the M.S. in Physics of Complex Systems in 2017 from Universit\`a  di Torino, Italy, and the PhD in Pure and Applied Mathematics in 2022 from Politecnico di Torino,Italy. He is currently a Research Assistant at the Department of  Mathematical Sciences, Politecnico di Torino, Italy. He was a Visiting Student at the Laboratory for Information and Decision Systems, Massachusetts Institute of Technology, in 2018--2020. His research focuses on control in network systems, with application to transportation and epidemics.
\end{IEEEbiography}

\begin{IEEEbiography}[{\includegraphics[width=1in,height=1.25in,clip,keepaspectratio]{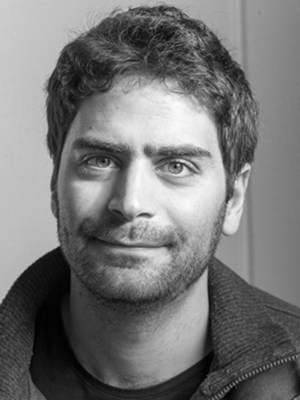}}]
{Giacomo Como}(M'12) is  a  Professor at  the Department  of  Mathematical  Sciences,  Politecnico di  Torino,  Italy. He is also a Senior Lecturer  at  the  Automatic  Control  Department, Lund  University,  Sweden.  He  received the B.Sc., M.S., and Ph.D.~degrees in Applied Mathematics  from  Politecnico  di  Torino, Italy, in  2002,  2004, and 2008, respectively. He was a Visiting Assistant in  Research  at  Yale  University  in  2006--2007  and  a Postdoctoral  Associate  at  the  Laboratory  for  Information  and  Decision  Systems,  Massachusetts  Institute of Technology in  2008--2011. Prof.~Como currently serves as Senior Editor for the \textit{IEEE Transactions on Control of Network Systems}, as Associate  Editor  for \textit{Automatica}, and as the  chair  of the  {IEEE-CSS  Technical  Committee  on  Networks  and  Communications}. He served as Associate Editor for the  \textit{IEEE Transactions on Network Science and Engineering} (2015-2021) and for the \textit{IEEE Transactions on Control of Network Systems} (2016-2022).  He was  the  IPC  chair  of  the  IFAC  Workshop  NecSys'15  and  a  semiplenary speaker  at  the  International  Symposium  MTNS'16.  
He  is  recipient  of  the 2015  George S. ~Axelby  Outstanding Paper Award.  His  research interests  are in  dynamics,  information,  and  control  in  network  systems. 
\end{IEEEbiography}

\begin{IEEEbiography}[{\includegraphics[width=1in,height=1.25in,clip,keepaspectratio]{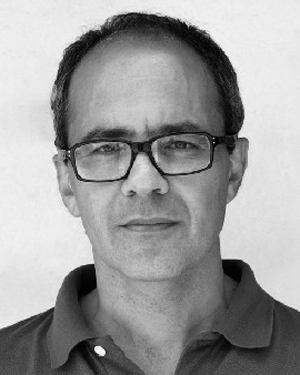}}]
{Fabio Fagnani}
received the Laurea degree in Mathematics from the University of Pisa and the Scuola Normale Superiore, Pisa, Italy, in 1986. He received the PhD degree in Mathematics from the University of Groningen,  Groningen,  The  Netherlands,  in 1991. From 1991 to 1998, he was an Assistant Professor at the Scuola Normale Superiore. In 1997, he was a Visiting Professor at the Massachusetts Institute of Technology. Since 1998, he has been with the Politecnico of Torino. From 2006 to 2012, he acted as Coordinator of the PhD program in Mathematics for Engineering Sciences at Politecnico di Torino. From 2012 to 2019, he served as the Head of the Department of Mathematical Sciences, Politecnico di Torino. His current research focuses on network systems, inferential distributed algorithms, and opinion dynamics. He is an Associate Editor of the \textit{IEEE Transactions on Automatic Control} and served in the same role for the \textit{IEEE Transactions on Network Science and Engineering} and the \textit{IEEE Transactions on Control of Network Systems}.
\end{IEEEbiography}
	
\end{document}